\documentclass[11pt]{article}

\usepackage[utf8]{inputenc}
\usepackage[english]{babel}
\usepackage{amsthm}
 \usepackage{amssymb}
 \usepackage{amsmath}
 \usepackage{hyperref}

\usepackage{comment}

\usepackage{graphicx}
\usepackage{color}
\usepackage{tikz}
\usepackage{float}
\usepackage{subfig}

\numberwithin{equation}{section}

\newtheorem{theorem}{Theorem}[section]
\newtheorem{definition}{Definition}[section]
\newtheorem{question}[theorem]{Question}
\newtheorem{lemma}[theorem]{Lemma}
\newtheorem{conjecture}[theorem]{Conjecture}
\newtheorem{proposition}[theorem]{Proposition}
\newtheorem{corollary}[theorem]{Corollary}
\newtheorem{remark}[theorem]{Remark}
\newtheorem{problem}[theorem]{Problem}

\def\bclaim{\begin{claim}}
\def\eclaim{\end{claim}}
\def\bdefin{\begin{definition}}
\def\edefin{\end{definition}}
\def\bcor{\begin{corollary}}
\def\ecor{\end{corollary}}
\def\bthm{\begin{theorem}}
\def\ethm{\end{theorem}}
\def\bconj{\begin{conjecture}}
\def\econj{\end{conjecture}}
\def\blem{\begin{lemma}}
\def\elem{\end{lemma}}
\def\blemma{\begin{lemma}}
\def\elemma{\end{lemma}}
\def\bprop{\begin{proposition}}
\def\eprop{\end{proposition}}
\def\bremark{\begin{remark}}
\def\eremark{\end{remark}}
\def\bprob{\begin{problem}}
\def\eprob{\end{problem}}
\def\bquestion{\begin{question}}
\def\equestion{\end{question}}

\newcommand{\PP}{{\mathbb P}} \newcommand{\RR}{\mathbb{R}}
 \newcommand{\CC}{{\mathbb C}}
\newcommand{\ZZ}{{\mathbb Z}} \newcommand{\NN}{{\mathbb N}}
\newcommand{\FF}{{\mathbb F}}

\def\a{\alpha} \def\be{\beta}
\def\o{\omega} 
\def\th{\theta} 
\def\vp{\varphi}

\newcommand{\dbar}{\bar\partial}
\newcommand{\ddbar}{\partial\dbar}

\def\Im{{\operatorname{Im}}}

\def\MA{Monge--Amp\`ere } 

\def\K{K\"ahler } \def\Kno{K\"ahler}
\def\KE{K\"ahler--Einstein } 
  \def\KEE{K\"ahler--Einstein edge }

\def\Ric{\hbox{\rm Ric}\,}

\def\h#1{\hbox{#1}}

\def\dis{\displaystyle}

\def\q{\quad} \def\qq{\qquad}

\def\ra{\rightarrow}
\def\pa{\partial}

\def\w{\wedge}
\def\i{\sqrt{-1}}

\def\sm{\setminus}
\def\al{\alpha}
\def\be{\beta}

\def\del{\partial}

\def\la{\lambda}
\def\beq{\begin{equation}}
\def\eeq{\end{equation}}
\def\beqno{\begin{equation*}}
\def\eeqno{\end{equation*}}
\def\eaeq{\end{aligned}}
\def\baeq{\begin{aligned}}

\def\bpf{\begin{proof}}
\def\epf{\end{proof}}

\def\eaeq{\end{aligned}}
\def\baeq{\begin{aligned}}

\chardef\inodot="10

\def\saldpno{strongly asymptotically log del Pezzo}

\def\oFS{\omega_{\FS}}

\def\la{\lambda}

\def\oCyl{\omega_{\Cyl}}

\def\er{\eqref}
\def\noi{\noindent}
\def\sms{\smallskip}

\def\be{\beta}
\def\K{K\"ahler }
\def\sm{\setminus}

\def\KEE{K\"ahler--Einstein edge }
\def\ra{\rightarrow}

\def\Bl{\operatorname{Bl}}
\def\FS{\operatorname{FS}}
\def\Cyl{\operatorname{Cyl}}

\def\lb#1{\label{#1}}
\def\t{\tau}

\title{Angle deformation of K\"ahler--Einstein edge metrics
\\
on Hirzebruch surfaces}

\date{30 June 2020}

\begin{document}
\author{Yanir A. Rubinstein, Kewei Zhang}
\maketitle

\centerline{\it Dedicated to Bernie Shiffman on the occassion of his retirement }

\begin{abstract}
We construct a family of \KEE metrics on all Hirzebruch surfaces using 
the Calabi ansatz 
and study their angle deformation. This allows us
to verify in some special cases a conjecture of Cheltsov--Rubinstein 
that predicts convergence towards a non-compact Calabi--Yau fibration
in the small angle limit.
We also give an example of a \KEE metric 
whose edge singularity is rigid, answering a question posed by Cheltsov.
\end{abstract}

\section{Introduction} 

Recently, Cheltsov--Rubinstein \cite{CR} put forward a 
conjectural picture in which non-compact Calabi--Yau fibrations
emerge as the small angle limit of families of
compact singular metrics known as \K edge metrics.
In a recent article, we verfied this picture in the most elementary setting of Riemann surfaces \cite{RZ}. Our goal in the present article is to give further evidence for the conjecture by verifying it for some special symmetric complex surfaces. In passing, we also answer a question of Cheltsov by proving the
existence of a rigid \KEE metric.

Let $M$ be a smooth complex projective manifold and let $D=D_1+\ldots+D_r$
be a simple normal crossing divisor in $M$.
Following Tian \cite[p. 147]{Tian1994}, a \K metric on $M$ with an edge singularity along a divisor $D\subset M$
is a smooth \K metric on $M\sm D$ that has a cone singularity
transverse to $D$ (`bent' at an angle $2\pi \be_i<2\pi$ along $D_i$). Such metrics are called
{\it \K edge metrics;}
we refer to \cite{JMR} for geometric analysis on such spaces
and to \cite{R14} for a detailled survey and further references.

In op. cit.  \cite{CR} one is concerned with small cone angle \K edge 
metrics that are also Einstein, called 
\KEE (KEE) metrics and of positive Ricci curvature.
For such metrics to exist it is necessary that
$(M,D)$ be {\it asymptotically log Fano}
\cite[Definition 1.1]{CR}, 
a positivity property generalizing
positivity of  the first Chern class (i.e., being Fano) in the sense
that $c_1(M)-\sum_{i=1}^r(1-\be_i(j))D_i$ is positive for a sequence 
$\be(j)=(\be_1(j),\ldots,\be_r(j))\in(0,1)^r$ tending to the origin
(being Fano corresponds to the other
`extreme' with $\be_1=\ldots=\be_r=1$ or alternatively there being no divisor, i.e., $D=0$). Assuming then
that $(M,D)$ is asymptotically log Fano and that it admits KEE metrics
$\{\o_{\be(j)}\}_{j\in\NN}$
with angles $\be_i(j)$ along $D_i$, it is conjectured there that by taking an appropriate limit these metrics will converge to a Calabi--Yau fibration
on the non-compact space $M\sm D$
 \cite[Conjecture 1.11]{CR}.

Conjecturally, the fibration may
 be described in terms of the 
 adjoint anticanonical linear system
 as follows. 
 Let $K_M$ denote the canonical bundle,
 so that $c_1(M)=-c_1(K_M)$ by definition.
 By our assumption $-K_M-D$ is nef as it is a limit of ample divisors.
 Let $d=\dim H^0(M,-K_M-D)$ and let
 $s_1,\ldots, s_d$ be a basis for the vector
 space $H^0(M,-K_M-D)$. The
 Kodaira map 
 \beq
 \lb{Kodairamapeq}
 M\ni z\mapsto [s_1(z):\ldots:s_d(z)]\in\PP^{d-1},
 \eeq
 is then a holomorphic map from $M$ 
 onto its image, a projective variety
 we denote by $Y\subset\PP^{d-1}$. 
 Denote
 by $\kappa:=
\dim Y$.
It was conjectured \cite[Conjecture 1.6]{CR} 
that $\kappa<\dim M$ (under the assumption that
the KEE metrics $\o_\be$ exist for small $\be$), or equivalently that
$(K_M+D)^{\dim M}=0$,
and this was established in dimension $2$
\cite{CR2}
and subsequently in all dimensions \cite{Fujita}.
The next step in the program described in \cite{CR} is to study the geometric limit, if such exists,
of such KEE metrics when the angles tend to zero.
More precisely, we are interested in the following conjecture \cite[Conjecture 1.8]{CR} 
(as customary, when discussing dimension 2 we will use `del Pezzo' instead of `Fano',
and replace $(M,D)$ by $(S,C)$): 
 
 \bconj
 \lb{CRconj}
 {\rm } 
Let $C=C_1+\ldots+C_r\not\sim-K_S$ be a disjoint union of smooth
curves in a projective surface $S$ such that $(S,C)$
is \saldpno. Let $\be\in(0,1)^r$ be a sequence tending
to the origin and suppose that there exist \KEE metrics
$\o_\be$ of angle $2\pi\be_i$ along $D_i$ and 
of positive Ricci curvature.
Then as $\be\ra0$ an appropriate
limit of $\o_\be$ converges to a
fibration of cylinders whose base is $\PP^1$.
 \econj

In fact, by general results of Kawamata and Shokurov $C$ can have at most two connected components, so that $r\in\{1,2\}$ \cite[Remark 3.7]{CR}. 
Note that part of the challenge is to identify what
exactly ``appropriate" means in Conjecture \ref{CRconj}. As we show
in this note the appropriate sense, at least in some cases, turns out to be 
rescaling the fibers of the fibration \er{Kodairamapeq}
in a precise way and considering a certain pointed Gromov--Haussdorff
limit, inspired by the asymptotic analysis in the Riemann surface case
\cite{RZ}. Moreover, we also show that the {\it un-rescaled} limit exists
as a collapsed Gromov--Haussdorff limit, and is the \KE metric on the base $Y$ (in this case the Fubini--Study
metric on $\PP^1$). We believe this should generalize to the setting of Conjecture
\ref{CRconj}.

In this note we concentrate on a particular family of $(S,C)$
as in Conjecture \ref{CRconj} and first construct
rather explicitly the sequence of KEE metrics $\o_\be$ for a sequence of small 
angles $\be\in(0,1)^2$ tending to zero.
Moreover, we analyze the small angle limit of this sequence and resolve 
Conjecture \ref{CRconj} for this family of pairs. The pairs we consider are
\beq
\lb{II2Aneq}
S=\FF_n,\, n\in\NN\cup\{0\},
\q 
C_1=Z_n,\, C_2\in|Z_n+nF|,
\eeq
where $\FF_n$ is the $n$-th Hirzebruch surface \cite{Hirzebruch}, $C_1$ is the $-n$-curve,
and $C_2$ is an $n$-curve; see \S\ref{compactifyingSec} and \S\ref{sec:kahler-class} for notation and precise definitions.
The `boundary' in these pairs consists of  
two disjoint components $C_1$ and $C_2$ which in some senses is an added difficulty, 
however, the pairs \er{II2Aneq} have the advantage of being toric  (by which we mean that 
$S$ is toric and each $C_i$ is a torus-invariant holomorphic curve in $S$), and in fact the pairs even admit the larger
symmetry group of the Calabi ansatz. 
The simplest sub-case, $n=0$, serves
as a general guide, and, in fact, as a nice bridge to our previous 1-dimensional
work \cite{RZ}. In the $n=0$ case, the boundary consists of two disjoint fibers in $\PP^1\times\PP^1$ and while the class
$-K_{\PP^1\times\PP^1}
-(1-\be_1)C_1-(1-\be_2)C_2$ is ample for all small $\be_1$
and $\be_2$, only the classes with $\be_1=\be_2$ admit KEE metrics which we denote
by $\o_{\be_1,\be_1}$. This is reminiscent of the situation for footballs, and
for a good reason. In fact the reason that only the classes with $\be_1=\be_2$
admit KEE metrics is that if $\be_1\not=\be_2$ then $\PP^1$ admits a conic
Ricci soliton football metric with angles $2\pi\be_1$ at the North pole and 
$2\pi\be_2$ at the South pole (see, e.g., \cite[\S2]{RZ}). Taking the product with the Fubini--Study metric
on $\PP^1$ we see that $(\PP^1\times\PP^1, C_1+C_2)$ admits a \Kno--Ricci soliton edge
metric with nontrivial vector field, hence cannot admit a KEE metric. 
When
$\be_1=\be_2$,
the KEE metric $\o_{\be_1,\be_1}$ 
is nothing but the product
of the volume one Fubini--Study metric on the base $Y=\PP^1$ and the 
volume $4\pi\be_1$ constant scalar curvature
conic metric, i.e., the football, on the fibers. Then, our previous
work \cite[Theorem 1.3]{RZ} shows that after rescaling the fiber metrics by $1/\be_1^2$
there is a pointed Gromov--Haussdorff limit that converges to
the product space $(\PP^1\times\CC^\star, \pi_1^*\oFS+C\pi_2^*\oCyl)$
for some $C>0$. In fact, this can also be shown by elementary complex
analysis (op. cit. is concerned with the more general
setting of Ricci solitons).

When $n>0$ the results of \cite{RZ} do not apply as the situation is no
longer a product one and we need different tools. 
First, we describe some of the well-known Calabi ansatz computations
that reduce the KEE equation to an ODE. 
These are classical when there are no edges, and as we show using the asymptotic
expansion of \cite[Theorem 1]{JMR} they generalize
naturally to the setting of edges.
Similarly to the case $n=0$
it turns out there is a relation between $\be_1$ and $\be_2$ that is
necessary and sufficient for the class 
$-K_{\FF_n}-(1-\be_1)C_1-(1-\be_2)C_2$ to admit a KEE representative
$\o_{\be_1,\be_2}$.
Moreover, $\be_2\ra0$ as $\be_1\ra0$ and moreover $\lim_{\be_1\ra0}\be_2/\be_1=1$. 
We then generalize the asymptotic analysis from $n=0$ to prove
Conjecture \ref{CRconj} for the pairs \er{II2Aneq} by analyzing rather explicitly
such sequences of KEE metrics $\o_{\be_1,\be_2}$ using the Calabi ansatz.

\bthm
\lb{mainthm}
Consider the pairs \er{II2Aneq}. Then
for each $n\in\NN$:

\sms
\noi
$\bullet\;$ There exists a family of KEE metrics $\o_{\be_1,\be_2}$ 
on the pairs \er{II2Aneq} for each 
$$
(\be_1,\be_2)
=
\Big(
\be_1,
\frac1{2n}\big(n\be_1-3+\sqrt{3(3-n\be_1)(1+n\be_1)}\big)
\Big),
\q \h{ for all $\be_1\in(0,2/n)\cap (0,1]$},
$$
with $\o_{\be_1,\be_2}$ cohomologous to
$$
\frac{2+n\be_2}{2-n\be_1}[C_2]-[C_1].
$$

\sms
\noi
$\bullet\;$
As $\be_1$ tends to zero, $(\FF_n,\o_{\be_1,\be_2})$ converges in 
the Gromov--Haussdorff sense to $(\PP^1,n\oFS)$,
where $\oFS$ is the Fubini--Study metric. 
Moreover, on the level of tensors,
as $\be_1$ tends to zero,
$\o_{\be_1,\be_2}$ restricted
to $\FF_n\setminus (C_1\cup C_2)$ converges
in all $C^k$ norms to a degenerate
tensor that is the pull-back
of $n\oFS$ on $\PP^1$ under 
the projection map to the 
zero section of the natural fibration $\pi_1:\FF_n\ra
\PP^1$.

\sms
\noi
$\bullet\;$
The fiberwise-rescaled metrics $\widetilde{\o_{\be_1,\be_2}}$
obtained by rescaling by $\be_1^{-2}$ only the restriction of 
$\o_{\be_1,\be_2}$ to each fiber (of the projection to $\PP^1$)
converges in $C^k$ on compact subsets as well as in the
pointed Gromov--Haussdorff sense to a cylinder fibration over $\PP^1$,
$(\PP^1\times\CC^\star,n\pi_1^*\oFS
+n\pi_2^*\oCyl)$, i.e.,
$$
\lim_{\be_1\ra0}\widetilde{\o_{\be_1,\be_2}}=
n(\pi_1^*\oFS
+
\pi_2^*\oCyl),
$$ 
where $\oCyl$ is the flat
metric on $\CC^\star$, and $\pi_i$ is the projection on to the $i$-th factor.
\ethm

\bremark
{\rm
In the third statement of Theorem
\ref{mainthm} there is a certain
choice of coordinates $(Z,w)$
on $\FF_n\setminus (C_1\cup C_2)$
(see \eqref{pi1pi2eq})
which then determines
projection maps $\pi_1$ and $\pi_2$. Thus, the limit
is only unique up to automorphisms
of $\PP^1\times\CC^*$.
}
\eremark

In 2015, Cheltsov posed the following question \cite{Cheltsov}.

\bquestion
\lb{VanyaQ}
Let $D$ be a smooth connected divisor in a smooth variety $X$,
and let $\be\in(0,1)$.
Does there exists a triple $(X,D,\be)$ 
such that $(X,D)$ admits a KEE metric of positive Ricci curvature
and angle $2\pi\be$ along
$D$ but does not admit KEE metrics of angles $2\pi\al$ along
$D$ for any $(0,1)\ni \al\not=\be$?
\equestion

Observe that Cheltsov's question
becomes interesting in dimension two and above since
there does not exist a constant scalar curvature conic Riemann sphere with
a single cone point. 
We answer question \ref{VanyaQ}  affirmatively.

\bcor
\lb{vanyathm}
Let $Z_{-1}$ be a smooth curve of self-intersection 1 disjoint from
the $-1$-curve in the first Hirzebruch surface $\FF_1$. 
The pair $(\FF_1,Z_{-1})$  admits a KEE metric of angle
$2\pi\be\in(0,2\pi]$ along $Z_{-1}$ if and only if
$\be=\sqrt3-1$.
\ecor

\bpf
The existence statement
for $n=1$ and $\be_1=1$
follows from
Theorem \ref{mainthm}
(or
from Corollary \ref{KEECor})
which
implies
that then $\be_2=\sqrt3-1$.
This gives a KEE metric with angle $2\pi(\sqrt3-1)$ along $Z_{-1}$
and smooth elsewhere.
This
KEE metric {\it cannot} be deformed to have an edge singularity
along the infinity section. 
Indeed, by \cite[Example 2.8]{CR2} the pair $(\FF_1,Z_{-1})$ does not admit KEE metrics
for {\it any angle smaller} than $\sqrt 3-1$ because $Z_{-1}$ log slope destabilizes it,
and neither for {\it any angle larger} than $\sqrt 3-1$ because $Z_1$ log slope destabilizes it. This proves Theorem \ref{vanyathm}.
\epf

The metric of Corollary 
\ref{vanyathm} is quite remarkable in that it is in fact perhaps the first example of an ``isolated" KEE
metric whose edge singularity cannot be deformed at all.

\subsection{Organization}

In Section \ref{sec:def-Hirz}, we provide several useful viewpoints of Hirzebruch surfaces for the reader's convenience. In Section \ref{sec:edge-metric}, we characterize \K edge metrics on Hirzebruch surfaces
(Proposition \ref{CalabiFnProp}) using the Calabi ansatz and the asymptotic analysis of \cite{JMR}. In Section \ref{sec:Einstein-eq}, we solve explicitly the K\"ahler--Einstein edge equation and determine the corresponding angles along the boundary divisors (Corollary \ref{KEECor}).
We emphasize that the reader that is only interested in Corollary \ref{vanyathm} can skip
Section \ref{sec:small-angle}
as Corollary \ref{vanyathm}
follows directly from Corllary \ref{KEECor}. 
Finally, in Section \ref{sec:small-angle}, we study the small angle limits of the K\"ahler--Einstein edge metrics
and prove Theorem \ref{mainthm}.

\section{Several descriptions of Hirzebruch surfaces}
\label{sec:def-Hirz}

\subsection{The quotient singularity}

Let $n\in\NN$. The simplest singularity in complex geometry is the quotient singularity constructed as follows.
Let $\ZZ_n$ act on $\CC^2$ by 
the
diagonal action, so
that the $\ZZ_n$-orbit of a point
$(a_1,a_2)$ is the 
collection of $n$ points
$$
\{(e^{2\pi \i l/n}a_1,e^{2\pi \i  l/n}a_2)
\in\CC^2
\,:\,
l=0,\ldots,n-1
\}.
$$
This is always a collection of $n$ distinct points, unless 
$(a_1,a_2)=(0,0)$.
So if we consider the orbit space (or
quotient space)
$$
\CC^2/\ZZ_n
$$
defined, as a set, as $\CC^2$ quotiented by the equivalence relation
\beq
\baeq
\label{Zkaction}
&(a_1,a_2)\sim (b_1,b_2)
\hbox{\ if there is
some $l\in\{0,\ldots,n-1\}$}\cr
&\hbox{
such that
$a_1=e^{2\pi \i l/n}b_1$ 
and 
$a_2=e^{2\pi \i l/n}b_2$ 
(with the same $l$),
}
\eaeq\eeq
we obtain an orbifold, smooth on the complement of a single singular point.

\subsection{Blow-up description of the total space}

In the remainder of this note we always assume
$$
n\in\NN,
$$
since the case $n=0$ was treated in the Introduction.
We resolve the quotient singularity
by blowing up $
\CC^2/\ZZ_n
$
at the single singular point.
Next we give an alternative global description of the resulting resolution.

Denote by 
\beq
\lb{blowupCnZk}
\Bl_0(\CC^2/\ZZ_n)
\eeq
the blow-up at the origin of (the orbifold, if $n\ge2$) $\CC^2/\ZZ_n$. Denote by 
$
H\ra\PP^{1}
$ the 
hyperplane bundle over $\PP^{1}$
and by 
$
H^{*}\ra\PP^{1}
$  the dual (tautological) line bundle. 
Then by
\beq
\lb{minusKspaceeq}
-nH_{\PP^{1}}
\eeq
we denote
the ($2$-dimensional) 
total space of $-nH\equiv (H^{*})^{\otimes n}$ 
considered as a line bundle over $\PP^{1}$. 
The following is an elementary exercise.
\blemma
\lb{blowupbihollemma}
$\Bl_0(\CC^2/\ZZ_n)$ is biholomorphic to
$
-nH_{\PP^{1}}.
$
\elemma

\subsection{Compactifying at infinity}
\lb{compactifyingSec}

By adding a point at infinity in each $\CC$ fiber (in the $-kH_{\PP^{1}}$ description) we obtain a compact space, a $\PP^1$ fibration over $\PP^{1}$.
The space is equivalently obtained by 
taking a product of each $\CC$ fiber
with another copy of $\CC$ and then 
taking the quotient under the $\CC^*$ action on the $\CC^2$ fibers (but
not acting on the base). This space is often
denoted by
\beq
\lb{Fnkdescreq}
\mathbb{P}
\big(
-nH_{\PP^{1}}\oplus\CC_{\PP^{1}}
\big),
\eeq
and we will simply denote it by
\beq
\lb{Fnkdefeq}
\mathbb{F}_{n},
\eeq
for any $n\ge0$.
In effect, we have added a copy of $\PP^{1}$. We therefore have two distinguished divisors in the space: the zero section 
\beq
\lb{Eeq}
Z_{n}:=\{w=0\} 
\eeq
(the exceptional divisor in the blow-up description) and the infinity section
\beq
\lb{Geq}
Z_{-n}:=\{w=\infty\}.
\eeq
By construction the two do not intersect. 
We emphasize that 
$$
\hbox{\it $Z_n$ is the $-n$-curve and $Z_{-n}$
is the ``section at infinity";}
$$
$Z_n$ can be contracted to give rise to the weighted projective space $\PP(1,1,n)$.

We can take \eqref{Fnkdescreq}
as the definition of 
$\mathbb{F}_{n}$. In that case
it actually makes sense for any
$n\in\ZZ$. And, since
for any vector bundle $A$ and line bundle
$L$ we have $\PP(A\otimes L)=\PP(A)$, it
follows that
$\mathbb{F}_{-n}$
is biholomorphic to $\mathbb{F}_{n}$ (take $L=2nH_{\PP^{1}}$)
with the biholomorphism exchanging the
zero and the infinity sections (i.e., $Z_n$ with $Z_{-n}$).
Recalling Lemma \ref{blowupbihollemma}
we have shown the following:

\blemma
$\mathbb{F}_{n}\setminus Z_{-n}$
is biholomorphic to 
$\Bl_0(\CC^2/\ZZ_n)$.
\elemma


\blemma 
\lb{betty2lem}
$\dim H^{1,1}_{\dbar}(\FF_{n})=2$.
\elemma 

\bpf
By \cite{Hirzebruch},
\begin{equation*}
\FF_n\xrightarrow[\operatorname{diffeo.}]{\cong}
    \begin{cases}
    S^2\times S^2,\ & n\text{ even},\\
    \PP^2\#\overline{\PP^2},\ & n\text{ odd}.\\
    \end{cases}
\end{equation*}
So one has $\dim H^2(\FF_n,\CC)=2$.
Thus
$
2\dim H^{2,0}_{\dbar}(\FF_{n})
+
\dim H^{1,1}_{\dbar}(\FF_{n})=2,
$
so it follows that $\dim H^{1,1}_{\dbar}(\FF_{n})$ is either 0
or 2. The former possibility is excluded since $Z_n\cong \PP^{1}$ is a non-trivial holomorphic $(1,1)$-cycle.
\epf

\bremark
\lb{formsinclassRem}
{\rm
We will explicitly construct (1,1)-forms in Section \ref{sec:kahler-class} below which gives an alternative direct description why 
$\dim H^{1,1}_{\dbar}(\FF_{n},\ZZ)\ge2$
hence equal to 2. In fact, we will construct explicit representatives for each K\"ahler class.
}
\eremark

\section{Edge metrics on Hirzebruch surfaces}
\label{sec:edge-metric}

\subsection{
Edge metrics on powers of the tautological line bundle
} 

In this section we elaborate on Remark
\ref{formsinclassRem} and give a geometric approach to 
compactification. 

We start by constructing some $U(2)$-invariant \K edge metrics on 
$-nH_{\PP^{1}}$ \eqref{minusKspaceeq}.
The construction goes back to Calabi \cite{Calabi,CalabiI} who
considered the case of smooth \K metrics.
The
generalization to the edge case is
not much harder.

Consider the Hermitian metric $h$
on $-H$ that assigns to each point in the total
space $(Z,w)$ the norm squared
$$
|(Z,w)|_{h}^2:=|w|^2||Z||^2=|w|^2(|Z_1|^2+|Z_2|^2),
$$
i.e., locally 
$$
h(Z)=||Z||^2,
$$
where $||Z||$ is the Euclidean norm of the vector $(Z_1,Z_2)$ in $\CC^2$.
Similarly, $h^n$ is a metric on $-nH$ and 
$$
|(Z,w)|_{h^n}^2:=|w|^2||Z||^2=|w|^2(|Z_1|^2+|Z_2|^2)^n,
$$
and 
$$
h^k(Z)=(|Z_1|^2+|Z_2|^2)^k.
$$ 
On the chart $Z_2\not=0$, we choose local 
holomorphic coordinates 
$z:=Z_1/Z_2,$
so that 
\beq
\lb{Zzeq}
Z=[Z_1:Z_2]=[z:1].
\eeq
The curvature of this metric is a \K form on the 
base (i.e.,
the zero section $Z_n\cong\PP^1$) which is given by
$$
\baeq
-\i\ddbar\log h^n
&=-n\i\ddbar\log||Z||^2
=-n\i\ddbar\
\log(1+|z|^2)^2
\cr
&=
-n(1+|z|^2)^{-2}
\i dz \w \overline{dz}
=
-n\oFS.
\eaeq
$$
In the sequel, we will use the above curvature property of
$h^n$, as well as the fact that it is a globally defined 
$U(2)$-invariant smooth function on $\Bl_0(\CC^2/\ZZ_n)$.

\subsection{The Calabi ansatz on the total space}

We use the logarithm of the global invariant function
from above as our coordinate from now on. That is,
we set
\beq
\lb{sdefeq}
s(Z,w):=
\log|(Z,w)|_{h^n}^2=\log|w|^2
+n\log(1+|z|^2)
, \q (Z,w)\in -nH_{\PP^1},
\eeq
and seek canonical \K metrics on $\Bl_0(\CC^2/\ZZ_n)$ that depend solely on $s$, namely, \K metrics of the form
$$
\eta=\sqrt{-1}\partial\bar{\partial}f(s),
$$
where $f$ is a smooth function. Our goal will be to 
determine appropriate $f$ that make $\eta$
have various desirable curvature properties and edge type
singularities along $Z_{\pm n}$.

Denote
$$
f'(s):=\frac {df}{ds}.
$$
Note,
$$
\frac {\del s}{\del w} = 
\frac 1 w,
\q
\frac {\del s}{\del Z_i} = 
n\frac {\overline{ Z_i}}{||Z||^2}, \q i=1,2.
$$
Working on the chart $Z_2=1$ (recall \er{Zzeq}),
$$
\baeq
&\begin{pmatrix}
\dis\frac
{\del^2 f}
{\del w \del \bar w}
&
\dis\frac
{\del^2 f}
{\del w \del \bar z}
\cr
\dis\frac
{\del^2 f}
{\del z \del \bar w}
&
\dis\frac
{\del^2 f}
{\del z \del \bar z}
\cr
\end{pmatrix}
=
\begin{pmatrix}
\dis\frac{f''}{|w|^2}
& 
\dis
\frac{nf''{z}}{w(|z|^2+1)}
\cr
\dis
\frac{nf''\overline{z}}{\overline{w}(|z|^2+1)}
&
\dis
\frac{
nf'+n^2f''|z|^2}{(|z|^2+1)^2}
\cr
\end{pmatrix}
\eaeq
$$
Set 
\beq
\lb{pi1pi2eq}
\pi_1(Z,w):=z,
\q
\pi_2(Z,w):=w,
\eeq
and denote by $\oFS$ the Fubini--Study metric
on $\PP^{1}$
and the flat \K form on the cylinder
$\CC^\star$ by
\beq
\lb{oCyleq}
\oCyl:=\frac{\sqrt{-1}dw\wedge d\bar{w}}{|w|^2}.
\eeq
Thus,
\beq
\baeq
\lb{etaeq}
\eta
&=
nf'\pi_1^*\oFS+f''\pi_2^*\oCyl
\cr
&\quad
+
n\frac{f''}{1+|z|^2}
\Big(
\i \frac {dw}{w} \w z\,\overline{dz}
+
\i \overline{z}dz 
\w \frac {\overline{dw}}{ 
\bar w}
\Big)
\cr
&\quad
+
n^2f''(s)(1+|z|^2)^{-2}
\overline{z}dz\w 
z\,\overline{dz}.
\eaeq
\eeq
From this computation we see that the two key
quantities are $f'$ and $f''$ (rather than $f$
itself). Both of these must be positive, i.e.,
$f$ must be an increasing convex function of $s$.
Inspired by this, consider a Legendre type
change of variables going back to Calabi \cite{Calabi} (cf. \cite{HwangS};
see \cite{RZ} for
a reference most closely following our notation):
\beq
\lb{tausubsteq}
\tau=\tau(s):=f^{\prime}(s),\q \varphi=\varphi(\tau)=\varphi(\tau(s)):=f^{\prime\prime}(s), \q \tau\in \Im f'=f'(\RR).
\eeq
Setting
\beq
\lb{alphaeq}
\alpha:=
n\frac{\overline{z}dz}
{1+|z|^2}
,
\eeq
and using \eqref{tausubsteq},
we may rewrite \eqref{etaeq} as 
\beq
\baeq
\lb{eta2eq}
\eta
&=
n\tau\pi_1^*\oFS
+
\vp\Big(
\pi_2^*\oCyl
+\i\a\w\overline{\a}
+\i\a\w \overline{dw/w}
+\i dw/w\w \overline{\a}
\Big)
.
\eaeq
\eeq

\subsection{\K classes}
\label{sec:kahler-class}

By construction, the normal bundle of $Z_{\pm n}$ is $\mathcal{O}_{\mathbb{P}^{1}}(\mp n)=\mp nH_{\PP^1}$, 
so
$$
Z_{\pm n}^2=\mp n, \q Z_n.Z_{-n}=0.
$$ 
It follows from this and Lemma \ref{betty2lem} that $Z_{n}$ and 
$Z_{-n}$ generate
the Picard group, and thus by the Nakai--Moishezon
criterion the \K classes are precisely represented by $-xZ_n+yZ_{-n}$ with
$y>x>0$ (this also follows directly from the Calabi ansatz).
The divisor class 
\beq
\lb{Feq}
F:=\frac{1}{n}\big(Z_{-n}-Z_n\big),
\eeq
has zero self intersection and intersects $Z_{\pm n}$ exactly
at one point, hence represents the fibers of the projection to the base $Z_n$.
Setting
\beq
\lb{boundaryeq}
C_1=Z_n, \q C_2=Z_{-n},
\eeq
and using \er{Feq} we recover the notation \er{II2Aneq}.
The canonical class can be determined as follows. Write
$K_{\mathbb{F}_{n}}=-xZ_n+yZ_{-n}.
$
By Riemann--Roch,
$$
-2=(K_{\mathbb{F}_{n}}+Z_{\pm n}). Z_{\pm n}
=
((1-x)Z_n+yZ_{-n}).Z_n=-n(1-x)
=
(Z_n+(1+y)Z_{-n}).Z_{-n}=n(1+y),
$$
so $x=1-2/n, \, y=-1-2/n$, and 
\beq
\lb{KFneq}
-K_{\mathbb{F}_{n}}\sim 
\Big(1-\frac2n\Big)Z_{n}
+
\Big(1+\frac2n\Big)Z_{-n}
\sim
2Z_n+(n+2)F.
\eeq

\subsection{The angle constraint as 
boundary data}

We rewrite the angle constraint in terms of $\tau$.
Note that
the domain of $\tau$ is 
$$
\Im f'
=
(\inf f', \sup f'),
$$ 
and
$\tau$ must be positive so $\inf f'\ge0$.

\blemma
\lb{inffprimelemma}
Suppose that $\eta$ restricts to a Riemannian metric on $E$. Then,
$\inf f'>0$.
\elemma

\bpf
If $\inf f'=0$, then by \eqref{etaeq}
and \eqref{Eeq}
the restriction of $\eta$ to
$E$ is identically zero, which means
the zero section is collapsed to a point,
a contradiction.
\epf

To simplify computations we will henceforth
assume we are in the situation of Lemma
\ref{inffprimelemma} and
rescale $\eta$, equivalently $f$, by
a positive constant so that
\beq
\lb{inffprimeeq}
\inf f'=1. 
\eeq
This is equivalent to rescaling the 
\K class of $\eta$.
Thus the only
contribution to the Poincar\'e--Lelong 
formula aside from $\oCyl$ will be from any 
vanishing of $\vp$ along $\{w=0\}$
and along $\{1/w=0\}$.

By construction, $Z_{\pm n}=\{s=\mp\infty\}$.
By our normalization above 
$\tau$
ranges in the domain $(1,T)$ with
\beq
\lb{Etaueq}
Z_n=\{\tau=1\} \q 
Z_{-n}=\{\tau=T\}.
\eeq

\blem
\lb{TFinite}
$T<\infty$.
\elem
\bpf
By assumption $\eta$ is a \K edge metric on the pair 
 \er{II2Aneq}. 
 Restricting $\eta$ to a fiber (i.e., say, to the
vertical section $\{z=0\}$) using 
\er{eta2eq}, we get
an $S^1$-invariant metric
\begin{equation}
\lb{surfrevoleq}
g=\frac{1}{2\varphi(\tau)}d\tau^2+2\varphi(\tau)d\theta^2,
\end{equation}
where $w=e^{s/2+\i\th}$ is a coordinate on $\{z=0\}$.
This follows in the same way as in the 1-dimensional setting
\cite[Lemma 2.1]{RZ}.
Note that here we implicitly used the fact that $\eta$ is $U(2)$-invariant,
hence its restriction to any fiber is $S^1$-invariant.
Since for any \K edge metric the 
volume of a complex submanifold
is finite (in this case it is a 
cohomological constant)
and since by \er{surfrevoleq}
the volume form on the fiber
is simply $d\tau\wedge d\theta$,
it follows from \er{Etaueq} that the volume of the fiber
is $2\pi(T-1)$, 
that is finite if and only
if $T<\infty$.
\epf

By our assumption \eqref{inffprimeeq} and Lemma \ref{TFinite}, 
$\tau$
ranges in a domain $(1,T)$ with $T<\infty$. Thus,
$$
\lim_{s\ra\pm\infty}\frac{d\tau}{ds}=0,
$$
i.e., using  \er{tausubsteq},
\beq
\lb{vpzeroseq}
\vp(1)=\vp(T)=0.
\eeq

Next, we rewrite the angle constraint at the edges in
terms of $\tau$.

\bprop
\lb{CalabiFnProp}
$\eta$ is a \K edge metric on 
the pair  \er{II2Aneq} with angle $2\pi\be_1$
along $C_1$ and  $2\pi\be_2$
along $C_2$ if and only if 
\beq
\lb{vpinitialdataeq}
\vp(1)=0, \q \vp'(1)=\be_1,
\q
\vp(T)=0, \q \vp'(T)=-\be_2.
\eeq
\eprop

\bpf
By \er{vpzeroseq} it remains to determine the derivatives of
$\vp$ at $1$ and $T$.

Suppose first that $\eta$ is a \K edge metric with angles
as stated. It follows from \cite[Theorem 1, Proposition 4.4]{JMR} 
that $f$ has complete asymptotic expansions both
near $w=0$ and $w=\infty$. Let us concentrate
on the former first. The leading term in that expansion is $|w|^{2\be_1}$
and using \er{sdefeq},
$$
\begin{aligned}
\varphi
&\sim C_1+C_2|w|^{2\beta_1}+
(C_3\sin\th+C_4\cos\th)|w|^2+O(|w|^{2+\epsilon})
\cr
&=
C_1+C_2e^{\beta_1s}+
(C_3\sin\th+C_4\cos\th)e^{s}+O(e^{(1+\epsilon)s})
\end{aligned}
$$
(note that $r$ in \cite[(56)]{JMR} is equal to 
$|w|^{\be_1}/\be_1$ in our notation,
see \cite[p. 102]{JMR}).
Note that $C_1=0$ by \er{vpzeroseq} (actually also $C_3=C_4=0$ as $\vp$
is independent of $\th$ but we do not need this).
Moreover, the expansion can be differentiated term-by-term 
as $|w|\rightarrow0$ or $s\ra-\infty$.
As $
\vp'(\tau)=\frac{\pa\vp}{\pa s}\frac{ds}{d\tau}
=\frac{\pa\vp}{\pa s}/\varphi,
$
we obtain
\begin{equation}\label{eq:0}
    \varphi(1)=0,\quad \varphi^\prime(1)=\beta_1.
\end{equation}
The same arguments imply that
\begin{equation}
    \label{eq:b1+b2}
    \varphi(T)=0,
    \quad \varphi^\prime(T)=-\beta_2,
\end{equation}
the minus sign coming from the fact that the leading term in the expansion is now 
$1/|w|^{2\be_2}=e^{-\be_2 s}$.

Conversely, suppose that \er{vpinitialdataeq} holds.
Then near $Z_n=\{\tau=1\}=\{w=0\}$ there exists a positive
smooth
function $F(z,w)$ and a positive constant $\delta$ so that
\beq
\lb{angle1eq}
\vp(\tau(z,w))=|w|^{2\delta}F(z,w).
\eeq
Using \eqref{sdefeq}, $|w|^{2\delta}=e^{\delta s}/(1+|z|^2)^{\delta}$,
and 
\beq
\lb{dtauds2eq}
\frac{d}{d\tau}=\frac1\vp\frac{d}{ds},
\eeq
so in the notation of \eqref{angle1eq}, 
$$
\baeq
\vp_\tau(1)
&=\lim_{s\ra-\infty}
\frac{
\frac{d}{ds}\Big(F(s)e^{\be s}\Big)
}
{
F(s)e^{\be s}
}
\cr
&=\lim_{s\ra-\infty}
\frac{
F'(s)e^{\be s}+\be F(s)e^{\be s} 
}
{
F(s)e^{\be s}
}
\cr
&=
\be+\lim_{s\ra-\infty} F'(s)/F(s).
\eaeq
$$
By assumption, $\lim_{s\ra-\infty} F(s)$
is some positive (finite) number, in particular 
we must have $\lim_{s\ra-\infty} F'(s)=0$
Thus, $d\vp/d\tau(1)=\delta$ and 
\beq
\lb{deltabeta}
\delta=\beta_1.
\eeq

Now, taking the top wedge product of \eqref{eta2eq}
gives,
\beq
\lb{topwedgeetaeq}
\baeq
\eta^2
&=
2\bigg(n\tau\vp(\tau)
\pi_1^*\oFS\w \pi_2^*\oCyl
+
\vp^2\Big(
\pi_2^*\oCyl\w\i\a\w\overline{\a}
+\i\a\w \overline{dw/w}\w\i dw/w\w \overline{\a}
\Big)\bigg)
\cr
&=
2n\tau\vp(\tau)
\pi_1^*\oFS\w \pi_2^*\oCyl
.
\eaeq
\eeq
Now \er{oCyleq}, \er{angle1eq} and \er{deltabeta} imply that
$\eta$ satisfies a complex \MA equation with right-hand side
equal to a smooth volume form times  $1/|w|^{2-2\be_1}$ near $Z_n$. Thus,
 \cite[Theorem 1]{JMR} applies and $\eta$ is a \K edge
metric with a complete asymptotic expansion near $Z_n$
(note that op. cit. is stated in the case
of a smooth connected divisor but applies verbatim in the
case of smooth disjoint divisors). The same
arguments apply near $Z_{-n}$ to conclude.
\epf

\section{The Einstein constraint}
\label{sec:Einstein-eq}

From the proof of Proposition \ref{topwedgeetaeq} we can derive
a formula for the Ricci tensor of $\eta$. Indeed, by \er{topwedgeetaeq}
and the Poincar\'e--Lelong formula, 
\beq
\lb{Riccomp}
\Ric\eta=
(1-\be_1)[C_1]+
(1-\be_2)[C_2]+
2\pi_1^*\oFS-\i\ddbar\log\tau
-\i\ddbar\log\vp
,
\eeq
where the last term is understood
to be the restriction of $\i\ddbar\log\vp$
to the complement of $C_1+C_2$.
For the remaining terms, compute
\beq
\baeq
\lb{Riccomp1}
\del_w\log f'
&=
\frac{f''}{f'}
\frac{\del s}{\del w}
=
\frac{f''}{wf'}
\cr
\del_{z}\log f'
&=
\frac{f''}{f'}
\frac{\del s}{\del z}
=
n
\frac{f''}{f'}
\frac{\overline{z}}{1+|z|^2}
\cr
\del_{w\bar w}\log f'
&=
\bigg(\frac{f''}{f'}\bigg)'
\frac 1w\frac{\del s}{\del \bar w}
=
\bigg(\frac{f''}{f'}\bigg)'\frac1{|w|^2}
\cr
\del_{z\bar w}\log f'
&=
n
\bigg(\frac{f''}{f'}\bigg)'
\frac{\overline{z/w}}{1+|z|^2}
\cr
\del_{z\bar z}\log f'
&=
n^2
\bigg(\frac{f''}{f'}\bigg)'
\frac{|z|^2}{(1+|z|^2)^2}
+n\frac{f''}{f'}
\frac{
1
}
{(1+|z|^2)^2}
,
\eaeq
\eeq
and this can be simplified by noting that
\beq
\lb{dtaudseq}
\frac{d \tau}{ds}=f''=\vp(\tau).
\eeq
Thus,
\beq
\baeq
\lb{Riccomp2}
\del_w\log f'
&=
\frac{\vp}{w\tau}
\cr
\del_{z_i}\log f'
&=
n
\frac{\vp}{\tau}
\frac{\overline{z}}{1+|z|^2}
\cr
\del_{w\bar w}\log f'
&=
\Big(\frac{\vp}{\tau}\Big)_\tau\frac\vp{|w|^2}
\cr
\del_{z\bar w}\log f'
&=
n
\Big(\frac{\vp}{\tau}\Big)_\tau\vp
\frac{\overline{z/w}}{1+|z|^2}
\cr
\del_{|z|^2}\log f'
&=
n^2
\Big(\frac{\vp}{\tau}\Big)_\tau\vp
\frac{|z|^2}{(1+|z|^2)^2}
+n\frac{\vp}{\tau}
\frac{
1}
{(1+|z|^2)^2}
,
\eaeq
\eeq
Next, $\i\ddbar\log\vp$ is computed similarly by
replacing $f'$ by $f''$ everywhere. A simplification
is obtained by noting that 
$$
\frac{f'''}{f''}=
\frac{\vp_\tau\vp}{\vp}=\vp_\tau.
$$
Thus,
\beq
\baeq
\lb{Riccomp3}
\del_w\log f''
&=
\frac{\vp_\tau}{w}
\cr
\del_{z}\log f''
&=
n
\frac{{\vp_\tau}\overline{z}}{1+|z|^2}
\cr
\del_{w\bar w}\log f''
&=
\frac{\vp_{\tau\tau}\vp}{|w|^2}
\cr
\del_{z\bar w}\log f''
&=
n
\frac{\vp_{\tau\tau}\vp\overline{z/w}}
{1+|z|^2}
\cr
\del_{z\bar z}\log f''
&=
n^2
\frac{\vp_{\tau\tau}\vp|z|^2}
{(1+|z|^2)^2}
+n{\vp_\tau}
\frac{
1}
{(1+|z|^2)^2}
.
\eaeq
\eeq
Altogether,
\beq
\baeq
\i\ddbar\log\tau
&=
n\vp/\tau\pi_1^*\oFS
\cr
&\q\q\q
+\vp(\vp/\tau)_{\tau}
\Big(
\pi_2^*\oCyl
+
\i\a\w\overline{\a}
\cr
&\q\q\q\q\q\q\q\q\q\q+\i\a\w \overline{dw/w}
+\i dw/w\w \overline{\a}.
\Big),
\eaeq
\eeq
and
\beq
\baeq
\i\ddbar\log\vp
&=
n\vp_\tau\pi_1^*\oFS
\cr
&\q\q\q
+\vp\vp_{\tau\tau}
\Big(
\pi_2^*\oCyl
+\i\a\w\overline{\a}
\cr
&\q\q\q\q\q\q\q\q\q\q\,
+\i\a\w \overline{dw/w}
+\i dw/w\w \overline{\a}.
\Big),
\eaeq
\eeq
So,
with $\alpha$ given by
\eqref{alphaeq}
\begin{equation}
\baeq
\label{expression of Ric}
\Ric\eta
&=
(1-\be_1)[C_1]+
(1-\be_2)[C_2]+
(2-n\vp/\tau-n\vp_\tau)\pi_1^*\oFS
\cr
&\q\q\q
-\vp(
\vp/\tau+
\vp_{\tau}
)_\tau
\Big(\pi_2^*\oCyl
+\i\a\w\overline{\a}
\cr
&\q\q\q\q\q\q\q\q\q\q\q\q\q\q\q
+\i\a\w \overline{dw/w}
+\i dw/w\w \overline{\a}
\Big).
\eaeq
\end{equation}
The Einstein edge equation 
\beq
\lb{KEEetaeq}
\Ric\eta=\la\eta+
(1-\be_1)[C_1]+
(1-\be_2)[C_2]
\eeq
becomes,
using \eqref{eta2eq} and \eqref{expression of Ric},
the pair of equations
$$
\baeq
2-n\vp/\tau-n\vp_\tau
&=n\la \tau
\cr
-\vp\big(
\vp/\tau+
\vp_{\tau}
\big)_\tau
&=\la \vp.
\eaeq
$$
Observe that the first equation implies the second
by differentiating in $\tau$. Also, setting $\tau=1$
and using \er{vpinitialdataeq}
implies
\beq
\lb{Lambdaeq}
\la=\frac 2n-\be_1.
\eeq
Observe that this already puts a constraint, as we must
require positive Ricci curvature on $\FF_n\sm C$, equivalently
\beq
\lb{be1constraint}
\be_1\in\Big(0,\frac2n\Big)\cap (0,1].
\eeq
Thus, the Einstein equation near $Z_n$ reduces to the
first-order initial value problem
\beq
\lb{Einsteintau2eq}
\vp_\tau+\frac\vp\tau=
\frac 2n + \Big(\be_1-\frac 2n\Big)\tau,
\qq \vp(1)=0.
\eeq

\subsection{Solving the Einstein equation}

Using an integration factor $\tau$ this becomes
\beq
\lb{Einsteintau2eq}
(\tau\vp)_\t
=
\frac 2n\tau
+ 
\Big(\be_1-\frac 2n\Big)\tau^2,
\qq \vp(1)=0,
\eeq
so
\beq
\lb{vptau1eq}
\vp
=
\frac 1{n}\frac{\tau^{2}-1}{\tau}
+ 
\frac1{3}\Big(\be_1-\frac 2n\Big)\frac{\tau^{3}-1}{\tau}
.
\eeq
Now, for this to correspond to a compact \K edge space with angle
$2\pi\be_2$ at $Z_{-n}$ we must satisfy \er{vpinitialdataeq}.
To that end, let us determine $T, \vp(T),$ and $\vp_\tau(T)$ in
\er{vpinitialdataeq} from \er{vptau1eq}. We factor $\vp$ as
\beq
\lb{vptaufactor}
\vp(\tau)=
\frac1{3}\Big(\be_1-\frac 2n\Big)
(\tau-1)(\tau-\al_1)(\tau-\al_2)/\tau,
\eeq
with $\al_1\le\al_2$.
Then,
$$
\frac1{3}\Big(\be_1-\frac 2n\Big)
(\tau-\al_1)(\tau-\al_2)=
\frac1{3}\Big(\be_1-\frac 2n\Big)(\tau^{2}+\tau+1)+
\frac 1{n}(\tau+1),
$$
so 
$$
-\al_1\al_2=\al_1+\al_2
=
\frac{1+n\be_1}{2-n\be_1}.
$$
By \er{be1constraint} we see that $\al_1<0<\al_2$,
so $T=\al_2$ if we can show $\al_2>1$. Solving the
quadratic equation for $\al_1,\al_2$ gives
\beq
\lb{al2eq}
\baeq
2\al_2
&=\frac{1+n\be_1+\big[(1+n\be_1)^2+4(1+n\be_1)(2-n\be_1)\big]^{1/2}}
{2-n\be_1}
\cr
&=\frac{1+n\be_1+[(9+6n\be_1-3n^2\be_1^2]^{1/2}}
{2-n\be_1}
.
\eaeq
\eeq
We claim that $\al_2>1$. According to \er{al2eq}, that amounts to verifying
$$
\frac{1+n\be_1+[(9+6n\be_1-3n^2\be_1^2]^{1/2}}
{2-n\be_1}
>2,
$$
which, after some manipulation precisely reduces to \er{be1constraint},
proving the claim.

Finally, it remains to compute 
the final equation in \er{vpinitialdataeq}, i.e., 
$\vp_\tau(\al_2)$. Using \er{vptaufactor}, 
\beq
\baeq
\lb{vptaufactorderiv}
\vp_\tau(\al_2)
&=
\frac1{3}\Big(\be_1-\frac 2n\Big)
(\tau-1)(\tau-\al_1)/\tau\Big|_{\tau=\al_2}
\cr
&=
\frac1{3}\Big(\be_1-\frac 2n\Big)
(\al_2-1)(\al_2-\al_1)/\al_2
.
\eaeq
\eeq
Note that
\beq
\lb{al2minusal1eq}
\baeq
{\al_2-1}
&=\frac
{\sqrt{3(3-n\be_1)(1+n\be_1)}+3(n\be_1-1)}
{4-2n\be_1},
\cr
\frac{\al_2-\al_1}
{\al_2}
&=\frac
{2\sqrt{3(3-n\be_1)(1+n\be_1)}}
{1+n\be_1+\sqrt{3(3-n\be_1)(1+n\be_1)}},
\eaeq
\eeq
Thus,
\beq
\lb{al3minusal1eq}
\baeq
\vp_\tau(\al_2)
&=
\frac{n\be_1-2}{3n}
\frac
{\sqrt{3(3-n\be_1)(1+n\be_1)}+3(n\be_1-1)}
{4-2n\be_1}
\frac
{2\sqrt{3(3-n\be_1)(1+n\be_1)}}
{1+n\be_1+\sqrt{3(3-n\be_1)(1+n\be_1)}}
\cr
&=
-\frac{1}{3n}
\frac
{
\sqrt{3(3-n\be_1)(1+n\be_1)}
\big(\sqrt{3(3-n\be_1)(1+n\be_1)}+3(n\be_1-1)\big)
}
{1+n\be_1+\sqrt{3(3-n\be_1)(1+n\be_1)}}
\cr
&=
-\frac{1}{3n}
\frac
{
3(3-n\be_1)(1+n\be_1)
+3(n\be_1-1)
\sqrt{3(3-n\be_1)(1+n\be_1)}
}
{1+n\be_1+\sqrt{3(3-n\be_1)(1+n\be_1)}}
\cr
&=
-\frac{1}{3n}
\bigg[
\frac
{
[3(3-n\be_1)(1+n\be_1)]^{3/2}
+3(n\be_1-1){3(3-n\be_1)(1+n\be_1)}
}
{-(1+n\be_1)^2+{3(3-n\be_1)(1+n\be_1)}}
\cr
&\qq\qq\qq
-\frac
{
3(3-n\be_1)(1+n\be_1)^2
+3(n^2\be_1^2-1)
\sqrt{3(3-n\be_1)(1+n\be_1)}
}
{-(1+n\be_1)^2+{3(3-n\be_1)(1+n\be_1)}}
\bigg]
\cr
&=
-\frac{1}{3n}
\bigg[
\frac
{
\sqrt{3(3-n\be_1)(1+n\be_1)}
\big[
3(3-n\be_1)(1+n\be_1)
-3(n^2\be_1^2-1)
\big]
}
{
-4n^2\be_1^2+8+4n\be_1
}
\cr
&\qq\qq\qq
+\frac{
3(3-n\be_1)(2n^2\be_1^2-4-2n\be_1)}
{
-4n^2\be_1^2+8+4n\be_1
}
\bigg]
\cr
&=
\frac{
3-n\be_1
}
{2n}
+
\frac{1}{n}
\frac
{
\sqrt{3(3-n\be_1)(1+n\be_1)}
\big[
(3-n\be_1)(1+n\be_1)
-(n^2\be_1^2-1)
\big]
}
{
4n^2\be_1^2-8-4n\be_1
}
\cr
&=
\frac{
3-n\be_1-\sqrt{3(3-n\be_1)(1+n\be_1)}
}
{2n}
=-\be_2.\eaeq
\eeq

Note that $\be_2=\be_2(\be_1)<\be_1$ for all $\be_1$.
Yet as $\be_1$ tends to zero, $\be_2$ tends to $\be_1$, to wit,
\beq
\baeq
\lb{be2asymp}
\be_2
&=
\frac{
n\be_1-3+3\sqrt{1+\frac23n\be_1-\frac13n^2\be_1^2}
}
{2n}
\cr
&=
\frac1{2n}
\bigg(
n\be_1-3+
3\Big(1+\frac13n\be_1
-\frac16n^2\be_1^2
-\frac18
\Big(
\frac23n\be_1-\frac13n^2\be_1^2
\Big)^2+O(\be_1^3)
\Big)
\bigg)
\cr
&=
\be_1-\frac n{3}\be_1^2+O(\be_1^3)
.
\eaeq
\eeq
Combining all the above together with
\er{KFneq}, \er{Lambdaeq}, and \er{KEEetaeq} we have shown
the following.

\bcor
\lb{KEECor}
Let $n\in\NN$.
For each $\be_1\in(0,2/n)\cap (0,1]$,
there exists a \KEE metric
$\o_{\be_1,\be_2}$
cohomologous to
$$
\frac{2+n\be_2}{2-n\be_1}[Z_{-n}]-[Z_n]
\sim
\frac{n(\be_1+\be_2)}{2-n\be_1}
[Z_{n}]
+
n\frac{2+n\be_2}{2-n\be_1}
[F]
$$
on the pair \er{II2Aneq}
with angles $2\pi\be_1$ along $C_1=Z_n$ and 
$2\pi\be_2=\pi\big(n\be_1-3+\sqrt{3(3-n\be_1)(1+n\be_1)}\big)/n$ along
$C_2=Z_{-n}$. 
One has $\be_2<\be_1$ and $\lim_{\be_1\ra0}\be_1/\be_2=1$. 
\ecor

\section{Small angle limits}
\label{sec:small-angle}

In this section 
we prove Theorem \ref{mainthm}.

First, let us determine the un-rescaled small angle limit of the KEE metrics
$\o_{\be_1,\be_2}$.
By \er{al2minusal1eq}
\beq
\lb{al4minusal1eq}
\baeq
{\al_2}
&=
1
+
3\frac
{\sqrt{1+\frac23n\be_1-\frac13n^2\be_1^2}+n\be_1-1}
{4-2n\be_1}
\cr
&=
1
+
3\frac
{
\frac43n\be_1
-\frac29n^2\be_1^2
+O(\be_1^3)
}
{4-2n\be_1}
=
1+n\be_1+\frac {n^2}3\be_1^2+O(\be_1^3),
\eaeq
\eeq
while,
\beq
\lb{al5minusal1eq}
\baeq
{\al_1}
&=
1
+
3\frac
{-\sqrt{1+\frac23n\be_1-\frac13n^2\be_1^2}+n\be_1-1}
{4-2n\be_1}
\cr
&=
1
+
3\frac
{
-2+\frac23n\be_1
+\frac29n^2\be_1^2
+O(\be_1^3)
}
{4-2n\be_1}
=
-\frac12+\frac n2\be_1+\frac {5n^2}{12}\be_1^2+O(\be_1^3)
.
\eaeq
\eeq
Recall that $\tau\in(1,\al_2)$. Thus, by \er{al4minusal1eq} 
we have $|\tau-1|=O(\be_1)$. Combining this with
\er{vptaufactor}, \er{al5minusal1eq}, 
and since 
by Corollary \ref{KEECor} and \er{eta2eq}, we have that
\beq
\lb{oKEEetaeq}
\o_{\be_1,\be_2}
=
n\tau\pi_1^*\oFS
+
\vp\Big(
\pi_2^*\oCyl
+\i\a\w\overline{\a}
+\i\a\w \overline{dw/w}
+\i dw/w\w \overline{\a}
\Big)
,
\eeq
we conclude that
\beq
\lim_{\be_1\ra0}
\o_{\be_1,\be_2}=
n\pi_1^*\oFS,
\eeq
and the convergence of tensors occurs smoothly,
so $(\FF_n\sm (Z_n\cup Z_{-n}),\o_{\be_1,\be_2})$ 
converges in the Gromov--Haussdorff sense to $(\PP^1,n\oFS)$.
(This combination of smooth convergence of tensors and simultaneous
collapse in the Gromov--Haussdorff sense is reminiscent of 
\cite[Theorem 2.4 (ii) (a)]{PR}.)
This concludes the proof of the first two statements
in Theorem \ref{mainthm}. 

We now turn to the last statement
in Theorem \ref{mainthm}. 
Inspired by \cite[Lemma 3.1]{RZ}, we change variable from $\tau\in(1,\al_2)$ to 
$$
y:=\frac{\tau-1-\frac{n\be_1}2}{n\be_1^2/2},
$$
with $y\in\Big(-\frac1{\be_1},\frac1{\be_1}+O(1)\Big)$
(recall \er{al4minusal1eq}), with $y=0$ roughly corresponding
to the mid-section between $Z_n$ and $Z_{-n}$.
Thus, by \er{vptaufactor}, \er{al4minusal1eq}, and \er{al5minusal1eq}, 
and since $|\tau-1|=O(\be_1)$,
\beq
\lb{vptaufactor2}
\vp(y)=
\frac{2-n\be_1}{2n}
\Big(\frac{n^2\be_1^2}4-\frac{n^2\be_1^4}{4}y^2\Big)
+O(\be_1^3), \q \Big(-\frac1{\be_1},\frac1{\be_1}+O(1)\Big).
\eeq

To determine a fiberwise-rescaled limit, define 
the fiberwise-rescaled metric, where the rescaling only occurs
for the terms that have a well-defined restriction to each fiber, i.e.,
$$
\widetilde{\o_{\be_1,\be_2}}
:=
n\tau\pi_1^*\oFS
+
\frac 1{\be_1^2}\vp
\pi_2^*\oCyl
+
\vp\Big(
\i\a\w\overline{\a}
+\i\a\w \overline{dw/w}
+\i dw/w\w \overline{\a}
\Big).
$$
with $\vp$ given by \er{vptaufactor2}.
As in the proof of Lemma \ref{TFinite}, the restriction of 
$\o_{\be_1,\be_2}$ to a fiber is given by
$$
\frac{1}{2\vp(\tau)}d\tau\otimes d\tau+ 2\vp(\tau)d\th\otimes d\th,
$$
so the 
restriction of 
$\widetilde{\o_{\be_1,\be_2}}$ to a fiber is given by
\beq
\lb{fiberwiserescaled}
\frac{n^2\be_1^2}{8\vp(y)}dy\otimes dy+ 2\frac{\vp(y)}{\be_1^2}d\th\otimes d\th,
\eeq
The Gromov--Haussdorff limit is then diffeomorphic to $\PP^1\times\CC^\star$.
Indeed,
using \er{vptaufactor2},
the length, with respect to $\omega_{\beta_1,\beta_2}$, of the path on each fiber between the intersection point of the fiber with $Z_n$ and its midpoint, the intersection of the fiber with the smooth section $\{y=0\}$ is
\beq
\int_1^{1+n\be_1/2}
\frac{d\tau}{\sqrt{\vp(\tau)}}.
\eeq 
To estimate this, recall
\er{vptaufactor}
and set $\xi:=\tau-1$,
\beq
\begin{aligned}
&\int_1^{1+n\be_1/2}
\frac{d\tau}{\sqrt{\vp(\tau)}}\\
        &=\int_{0}^{\frac{n}{2}\beta_1}\frac{\sqrt{\xi+1}}{\sqrt{\frac{1}{3}(\beta_1-\frac{2}{n})\xi(\xi+1-\alpha_1)(\xi+1-\alpha_2)}}\;\textrm{d}\xi\\
        =&\int_{0}^{\frac{n}{2}\beta_1}\frac{\sqrt{\xi+1}}{\sqrt{\frac{1}{3}(\beta_1-\frac{2}{n})\xi(\xi+\frac{3}{2}-\frac{n}{2}\beta_1-\frac{5n^2}{12}\beta_1^2+o(\beta_1^2))(\xi-n\beta_1-\frac{n^2}{3}\beta_1^2+o(\beta_1^2))}}\;\textrm{d}\xi
        \label{eq: estimate dist}
\end{aligned}
\eeq
    As $\beta_1\to 0$, the term $\displaystyle \frac{\sqrt{\xi+1}}{\sqrt{\frac{1}{3}(\frac{2}{n}-\beta_1)}}\cdot \frac{1}{\sqrt{\xi+\frac{3}{2}-\frac{n}{2}\beta_1-\frac{5n^2}{12}\beta_1^2+o(\beta_1^2)}}$ in the integral is uniformly bounded. Thus, to estimate the last displayed equation we consider 
    $$
    \begin{aligned}
        &\int_{0}^{\frac{n}{2}\beta_1} \frac{\textrm{d}\xi}{\sqrt{-\xi(\xi-n\beta_1+o(\beta_1))}}\\
        \overset{u=\frac{\xi}{\beta_1}}{=\joinrel=}&\int_0^{\frac{n}{2}}\frac{\beta_1\;\textrm{d}u}{\sqrt{-\beta_1 u(\beta_1 u-\beta_1 n+o(\beta_1))}}\\
        \overset{\beta_1\to 0}{=\joinrel=}&\int_0^{\frac{n}{2}}\frac{\textrm{d}u}{\sqrt{u(n-u)}}=O(1).
    \end{aligned}
    $$
    Similarly we also get $O(1)$ for the distance between $Z_{-n}$ and $\{y=0\}$ with respect to $\omega_{\beta_1,\beta_2}$. Hence, after the rescaling the fiberwise metric by $\beta_1^{-2}$, these same distances must be $O(\beta_1^{-1})$, and 
    in the limit $\beta_1\ra0$ we must get the product differential structure on $\mathbb{P}^1\times \mathbb{C}^*$. 
Moreover, in the limit $\beta_1\to 0$, \er{fiberwiserescaled}
 converges pointwise on compact subsets to 
    \begin{equation*}
        \frac{n}{2}{dy\otimes dy}+\frac{n}{2}{d\theta\otimes d\theta}.
    \end{equation*}
    The limiting metric in the pointed Gromov--Hausdorff sense is
    \begin{equation*}
        \lim_{\beta_1\to 0}\widetilde{\omega_{\beta_1,\beta_2}}=n(\pi_1^*\omega_{\operatorname{FS}}+\pi_2^*\omega_{\operatorname{cyl}}).
    \end{equation*}%
%
%
whose \K form is $nd\zeta\w\overline{d\zeta}/|\zeta|^2=n\oCyl$ (recall \er{oCyleq}) with
$\zeta:=e^{y+\i\th}$.
Combining this with \er{oKEEetaeq}, \er{vptaufactor2} and \er{fiberwiserescaled}, the limiting metric, in the pointed
Gromov--Hausdorff sense is then 
$$
\lim_{\be_1\ra0}\widetilde{\o_{\be_1,\be_2}}=
n(\pi_1^*\oFS
+
\pi_2^*\oCyl)
.
$$
This concludes the proof of Theorem \ref{mainthm}.

\bremark
\lb{VanyaRem}
{\rm
We mention in passing an interesting borderline phenomenon that occurs
in the cases $n\in\{1,2\}$. By Corollary \ref{KEECor}, when $n=2$ the metrics $\o_{\be_1,\be_2}$ exists for {\it all} $\be_1\in(0,1)$.
Naturally,
one may ask about the {\it large angle limit} $\be_1\ra1$. It turns out that in this 
case $\al_2\ra\infty$ and one obtains Gromov--Haussdorff convergence to the complete Ricci flat Eguchi--Hanson 
metric on the non-compact space $-2H_{\PP^1}$ as $C_2=Z_{-n}$ gets pushed-off to infinity.
Similarly, when $n=1$, the metrics $\o_{\be_1,\be_2}$ exists for {\it all} $\be_1\in(0,2)$. In the limit $\be_1\ra2$ one has $\be_2\ra1$
and the limit is now a complete Ricci flat metric on $-H_{\PP^1}$ with angle $4\pi$ 
along $Z_{1}$. We discuss these examples in detail elsewhere \cite{RZ2}.
}
\eremark

\paragraph{Acknowledgments.} 

Thanks to I.A. Cheltsov for posing to us Question \ref{VanyaQ}
and to Y. Ji for a careful reading and several corrections.
Research supported by NSF grants DMS-1515703,1906370,
a UMD--FAPESP Seed Grant, 
the China Scholarship Council award 201706010020,
the China post-doctoral grant BX20190014
and the Rosi \& Max Varon Visiting
		Professorship at the Weizmann Institute of Science in Fall 2019 and Spring 2020 to which Y.A.R. is
		grateful for the excellent research conditions. K.Z. 
		was a Visiting Scholar at the University of Maryland in
		2017--2018 when this work was initiated and is grateful 
		for its excellent research conditions.

\def\bi#1{\bibitem{#1}}

\let\omegaLDthebibliography\thebibliography 
\renewcommand\thebibliography[1]{
  \omegaLDthebibliography{#1}
  \setlength{\parskip}{1pt}
  \setlength{\itemsep}{1pt plus 0.3ex}
}

\vbox{
\vspace{0.1in}
 \noindent {\sc University of Maryland}\\
 {\tt yanir@alum.mit.edu}

\bigskip
 \noindent {\sc Beijing International Center for Mathematical Research, Peking University}\\
 {\tt kwzhang@pku.edu.cn}
}

\end{document}